\vsize=9.0in\voffset=1cm
\looseness=2


\message{fonts,}

\font\tenrm=cmr10
\font\ninerm=cmr9
\font\eightrm=cmr8
\font\teni=cmmi10
\font\ninei=cmmi9
\font\eighti=cmmi8
\font\ninesy=cmsy9
\font\tensy=cmsy10
\font\eightsy=cmsy8
\font\tenbf=cmbx10
\font\ninebf=cmbx9
\font\tentt=cmtt10
\font\ninett=cmtt9

\font\ninesl=cmsl9
\font\eightsl=cmsl8

\font\nineit=cmti9
\font\eightit=cmti8

\skewchar\ninei='177 \skewchar\eighti='177
\skewchar\ninesy='60 \skewchar\eightsy='60

\def\eightpoint{\def\rm{\fam0\eightrm} 
\normalbaselineskip=9pt
\normallineskiplimit=-1pt
\normallineskip=0pt

\textfont0=\eightrm \scriptfont0=\sevenrm \scriptscriptfont0=\fiverm
\textfont1=\ninei \scriptfont1=\seveni \scriptscriptfont1=\fivei
\textfont2=\ninesy \scriptfont2=\sevensy \scriptscriptfont2=\fivesy
\textfont3=\tenex \scriptfont3=\tenex \scriptscriptfont3=\tenex
\textfont\itfam=\eightit  \def\it{\fam\itfam\eightit} 
\textfont\slfam=\eightsl \def\sl{\fam\slfam\eightsl} 

\setbox\strutbox=\hbox{\vrule height6pt depth2pt width0pt}%
\normalbaselines \rm}

\def\ninepoint{\def\rm{\fam0\ninerm} 
\textfont0=\ninerm \scriptfont0=\sevenrm \scriptscriptfont0=\fiverm
\textfont1=\ninei \scriptfont1=\seveni \scriptscriptfont1=\fivei
\textfont2=\ninesy \scriptfont2=\sevensy \scriptscriptfont2=\fivesy
\textfont3=\tenex \scriptfont3=\tenex \scriptscriptfont3=\tenex
\textfont\itfam=\nineit  \def\it{\fam\itfam\nineit} 
\textfont\slfam=\ninesl \def\sl{\fam\slfam\ninesl} 
\textfont\bffam=\ninebf \scriptfont\bffam=\sevenbf
\scriptscriptfont\bffam=\fivebf \def\bf{\fam\bffam\ninebf} 
\textfont\ttfam=\ninett \def\tt{\fam\ttfam\ninett} 

\normalbaselineskip=11pt
\setbox\strutbox=\hbox{\vrule height8pt depth3pt width0pt}%
\let \smc=\sevenrm \let\big=\ninebig \normalbaselines
\parindent=1em
\rm}

\def\tenpoint{\def\rm{\fam0\tenrm} 
\textfont0=\tenrm \scriptfont0=\ninerm \scriptscriptfont0=\fiverm
\textfont1=\teni \scriptfont1=\seveni \scriptscriptfont1=\fivei
\textfont2=\tensy \scriptfont2=\sevensy \scriptscriptfont2=\fivesy
\textfont3=\tenex \scriptfont3=\tenex \scriptscriptfont3=\tenex
\textfont\itfam=\nineit  \def\it{\fam\itfam\nineit} 
\textfont\slfam=\ninesl \def\sl{\fam\slfam\ninesl} 
\textfont\bffam=\ninebf \scriptfont\bffam=\sevenbf
\scriptscriptfont\bffam=\fivebf \def\bf{\fam\bffam\tenbf} 
\textfont\ttfam=\tentt \def\tt{\fam\ttfam\tentt} 

\normalbaselineskip=11pt
\setbox\strutbox=\hbox{\vrule height8pt depth3pt width0pt}%
\let \smc=\sevenrm \let\big=\ninebig \normalbaselines
\parindent=1em
\rm}

\message{fin format jgr}

\vskip 4 mm
\magnification=1200
\font\Bbb=msbm10
\def\R{\hbox{\Bbb R}}
\def\C{\hbox{\Bbb C}}
\def\N{\hbox{\Bbb N}}
\def\pa{\partial}
\def\v{\varphi}
\def\ep{\varepsilon}

\vskip 4 mm

\centerline{\bf Absence of exponentially localized solitons}
\centerline{\bf for the Novikov-Veselov equation at positive energy}

\centerline{\bf R.G. Novikov}
\vskip 4 mm
\noindent
{\ninerm CNRS (UMR 7641), Centre de Math\'ematques Appliqu\'ees, Ecole
Polytechnique,}

\noindent
{\ninerm 91128 Palaiseau, France, and}

\noindent
{\ninerm IIEPT RAS - MITPAN, Profsoyuznaya str., 84/32,  Moscow 117997,
Russia}

\noindent
{\ninerm e-mail: novikov@cmap.polytechnique.fr}

\vskip 4 mm
{\bf Abstract.}\ In this note we show that the
Novikov-Veselov equation (NV-equation) at positive energy
(an analog of KdV in
2+1 dimensions) has no exponentially localized solitons in the two-
dimensional sense.

\vskip 2 mm
{\it 1.Introduction and Theorem 1.}
We consider the following 2+1 - dimensional analog of the KdV equaion (
Novikov-Veselov equation):
$$\eqalign{
&\pa_tv=4Re\,(4\pa_z^3v+\pa_z(vw)-E\pa_zw),\cr
&\pa_{\bar z}w=-3\pa_zv,\ \ v=\bar v,\ \ E\in\R,\cr
&v=v(x,t),\ \ w=w(x,t),\ \ x=(x_1,x_2)\in\R^2,\ \ t\in\R,\cr}\eqno(1)$$
where
$$\pa_t={\pa\over \pa t},\ \
\pa_z={1\over 2}\bigl({\pa\over \pa x_1}-i{\pa\over \pa x_2}\bigr),\ \
\pa_{\bar z}={1\over 2}\bigl({\pa\over \pa x_1}+i{\pa\over \pa x_2}\bigr).
\eqno(2)$$
We assume that
$$\eqalign{
&v\ \ {\rm is\ sufficiently\ regular\ and\ has\
sufficient\ decay\ as}\ \ |x|\to\infty,\cr
&w\ \ {\rm is\ decaying\ as}\ \ |x|\to\infty.\cr}\eqno(3)$$
Equation (1) is contained implicitly in the paper of S.V.Manakov [M] as
an equation possessing the following representation
$${\pa (L-E)\over \pa t}=[L-E,A]+B(L-E)\eqno(4)$$
(Manakov   L-A-B- triple),
where $L=-\Delta+v(x,t)$, $\Delta=4\pa_z\pa_{\bar z}$, $A$ and $B$ are
suitable differential operators of the third and zero order respectively,
$[\cdot,\cdot]$ denotes the commutator. Equation (1) was written in an
explicit form by S.P.Novikov and A.P.Veselov in [NV1], [NV2], where
higher analogs of (1) were also constructed.
Note also that the both Kadomtsev-Petviashvili (KP) equations can be obtained
from (1) by considering an appropriate limit $E\to \pm\infty$, see [ZS], [G2].

For the case when
$$v(x_1,x_2,t),\  w(x_1,x_2,t)\ \ {\rm are\ independent\ of}\ \ x_2\eqno(5)$$
equation (1) is reduced to
$$\pa_tv=2\pa_x^3v-12v\pa_xv+6E\pa_xv,\ \ x\in\R,\ \ t\in\R.\eqno(6)$$
In terms of $u(x,t)$ such that
$$v(x,t)=u(-2t,x+6Et),\ \ x\in\R,\ \ t\in\R,\eqno(7)$$
equation (6) takes the standard form of the KdV equation (see [NMPZ]):
$$\pa_tu-6u\pa_xu+\pa_x^3u=0,\ \ x\in\R,\ \ t\in\R.\eqno(8)$$
It is well-known (see [NMPZ]) that (8) has the soliton solutions
$$u(x,t)=u_{\kappa,\v}(x-4\kappa^2t)=
-{2\kappa^2\over ch^2(\kappa(x-4\kappa^2t-\v))},\ \
x\in\R,\ \ t\in\R,\ \kappa\in ]0,+\infty [,\ \v\in\R.\eqno(9)$$

In addition, one can see that
$$\eqalign{
&u_{\kappa,\v}\in C^{\infty}(\R),\cr
&\pa_x^ju_{\kappa,\v}(x)=O(e^{-2\kappa |x|})\ \ {\rm as}\ \ x\to\infty,\
j=0,1,2,3,\ldots\cr}\eqno(10)$$

Properties (10) show, in particular, that the solitons of (9) are
exponentially localized in $x$.

In the present note we obtain, in particular,  the following result:
\vskip 2 mm
{\bf Theorem 1.}
{\sl
Let $v,w$ satisfy (1) for $E=E_{fix}>0$, where
$$\eqalign{
&v(x,t)=V(x-ct),\ \ x\in\R^2,\ \ c=(c_1,c_2)\in\R^2,\cr
&V\in C^3(\R^2),\ \ \pa_x^jV(x)=O(e^{-\alpha |x|})\ \ for |x|\to\infty,\ \
|j|\le 3\ \ and\ some\  \alpha>0,\cr}\eqno(11a)$$
(where $j=(j_1,j_2)\in (0\cup \N)^2$, $|j|=|j_1|+|j_2|$,
$\pa_x^j=\pa^{j_1+j_2}\big/\pa x_1^{j_1}\pa x_2^{j_2}$),
$$w(\cdot,t)\in C(\R^2),\ \ w(x,t)\to 0\ \ as\ \ |x|\to\infty,\ \ t\in\R.
\eqno(11b)$$

Then $V\equiv 0$, $v\equiv 0$, $w\equiv 0$.
}

Theorem 1 shows that equation (1)  for $E>0$ has no nonzero solitons
(travel wave solutions) exponentially
localized in $x$ in the two-dimensional sense.
For $E<0$ this result will be given in [KN]. Note also that some other
integrable systems in 2+1 dimensions admit exponentially decaying solitons
in all directions on the plane, see [BLMP], [FS].

The proof of Theorem 1 is based on  Proposition 1  and Proposition 2, see
Section 4.
In turn, Proposition 2 is based, in particular, on Lemma 1 and Lemma 2.

Lemma 1, Lemma 2 and Proposition 1 are recalled in Section 2. Proposition 2
is given in Section 3. It seems that the result of Proposition 2 (that
sufficiently localized travel wave solutions for the NV-equation (1) for
$E=E_{fix}>0$ have zero scattering amplitude for the two-dimensional
Schr\"odinger equation (12)) was not yet formulated in the literature.

\vskip 2 mm
{\it 2. Lemma 1, Lemma 2 and Proposition 1.}
Consider the equation
$$-\Delta\psi+v(x)\psi=E\psi,\ \ x\in\R^2,\ \ E=E_{fix}>0,\eqno(12)$$
where
$$\eqalign{
&v(x)=\overline{v(x)},\ \ x\in\R^2,\cr
&(1+|x|)^{2+\ep}v(x)\in L^{\infty}(\R^2)\ \ {\rm (as\ a\ function\ of}\ \
x\in\R^2)\ \ {\rm for\ some}\ \ \ep>0.\cr}\eqno(13)$$

It is known that for any $k\in\R^2$, such that    $k^2=E$, there exists
an unique continuous solution $\psi^+(x,k)$ of equation (12) with the
following asymptotics:
$$\psi^+(x,k)=e^{ikx}-i\pi\sqrt{2\pi}e^{-i\pi/4}f(k,|k|{x\over |x|})
{e^{i|k| |x|}\over \sqrt{|k| |x|}}+o({1\over \sqrt{|x|}})\ \ {\rm as}\ \
|x|\to\infty.\eqno(14)$$
This solution describes scattering of incident plane wave $e^{ikx}$ on the
potential $v$. The function  $f$ on
$${\cal M}_E=\{k\in\R^2,\ l\in\R^2:\ k^2=l^2=E\}\eqno(15)$$
arising in (14) is the scattering amplitude for $v$ in the framework of
equation (12).
Under assumptions (13), it is known, in particular, that
$$f\in C({\cal M}_E).\eqno(16)$$

\vskip 2 mm
{\bf Lemma 1.}
{\sl
Let $v$ satisfy (13) and $v_y$, $y\in\R^2$, be defined by
$$v_y(x)=v(x-y),\ \ x\in\R^2.\eqno(17)$$
Then the scattering amplitude $f$ for $v$ and the scattering amplitude
$f_y$ for $v_y$ are related by the formula
$$f_y(k,l)=f(k,l)e^{iy(k-l)},\ \ (k,l)\in {\cal M}_E,\ \
y=(y_1,y_2)\in\R^2.\eqno(18)$$
}

Lemma 1 follows, for example, from the definition of the scattering amplitude
by means of (14) and the fact that $\psi^+(x-y,k)$ solves (12) for $v$
replaced by $v_y$, where $k^2=E$.

Lemma 1 was given, for example, in [N3].

\vskip 2 mm
{\bf Lemma 2.}
{\sl
Let $v,w$ satisfy (1), (3), where $E=E_{fix}>0$. Then the scattering amplitude
 $f(\cdot,\cdot,t)$ for $v(\cdot,t)$ and the scattering amplitude
$f(\cdot,\cdot,0)$ for $v(\cdot,0)$ are related by
$$f(k,l,t)=f(k,l,0)\exp[2it(k_1^3-3k_1k_2^2-l_1^3+3l_1l_2^2)],\ \
(k,l)\in {\cal M}_E,\ \ t\in\R.\eqno(19)$$
}

Lemma 2 was given for the first time in [N1].

Note that in the framework of Lemma 2 properties (3) can be specified as
follows:
$$\eqalign{
&v,w\in C(\R^2\times\R)\ \ {\rm and\ for\ each}\ \ t\in\R\ \ {\rm the\
following\ properties\ are\ fulfiled:}\cr
&v(\cdot,t)\in C^3(\R^2),\ \ \pa_x^jv(x,t)=O(|x|^{-2-\ep})\ \ {\rm for}\ \
|x|\to\infty,\ \ |j|\le 3\ \ {\rm and\ some}\ \ \ep>0,\cr
&w(x,t)\to 0\ \ {\rm for}\ \ |x|\to\infty.\cr}\eqno(20)$$

\vskip 2 mm
{\bf Proosition 1.}
{\sl
Let
$$v(x)=\overline{v(x)},\ \ e^{\alpha |x|}v(x)\in L^{\infty}(\R^2)\ \ ({\rm
as\ a\ function\ of}\ \ x)\ \ {\rm for\ some}\ \ \alpha>0 \eqno(21)$$
and the scattering amplitude $f\equiv 0$ on ${\cal M}_E$ for this potential
for some $E=E_{fix}>0$. Then $v\equiv 0$ in $L^{\infty}(\R^2)$.
}

In the general case the result of Proposition 1 was given for the first time
in [GN].
Under the additional assumption that $v$ is sufficiently small (in comparison
with $E$) the result of Proposition 1  was given for the first time in
[N2]-[N4].

\vskip 2 mm
{\it 3. Transparency of solitons.}
In this section we show that sufficiently localized solitons (travel wave
solutions) for the NV-equation (1) for $E=E_{fix}>0$ have zero scattering
amplitude for the two-dimensional Schr\"odinger equation (12).

\vskip 2 mm
{\bf Proposition 2.}
{\sl
Let $v,w$ satisfy (1) for $E=E_{fix}>0$, where
$$\eqalign{
&v(x,t)=V(x-ct),\ \ x\in\R^2,\ \ c=(c_1,c_2)\in\R^2,\cr
&V\in C^3(\R^2),\ \ \pa_x^jV(x)=O(|x|^{-2-\ep})\ \ for |x|\to\infty,\ \
|j|\le 3\ \ and\ some\  \ep>0,\cr}\eqno(22a)$$
$$w(\cdot,t)\in C(\R^2),\ \ w(x,t)\to 0\ \ as\ \ |x|\to\infty,\ \ t\in\R.
\eqno(22b)$$

Then
$$f\equiv 0\ \ {\rm on}\ \ {\cal M}_E,\eqno(23)$$
where $f$ is the scattering amplitude for $v(x)=V(x)$ in the framework of
the Schr\"odinger equation (12).
}

The proof of Proposition 2 consists in the following.

We consider
$$T=\{\lambda\in\C:\ \ |\lambda|=1\}.\eqno(24)$$
We use that
$${\cal M}_E\approx T\times T,\ \ E=E_{fix}>0,\eqno(25)$$
where diffeomorphism (25) is given by the formulas:
$$\lambda={{k_1+ik_2}\over \sqrt{E}},\ \
\lambda^{\prime}={{l_1+il_2}\over \sqrt{E}},\ \ (k,l)\in {\cal M}_E,\eqno(26)
$$
$$\eqalign{
&k_1={\sqrt{E}\over 2}\bigl(\lambda+{1\over \lambda}\bigr),\ \
k_2={i\sqrt{E}\over 2}\bigl({1\over \lambda}-\lambda\bigr),\cr
&l_1={\sqrt{E}\over 2}\bigl(\lambda^{\prime}+{1\over \lambda^{\prime}}\bigr),
\ \
l_2={i\sqrt{E}\over 2}\bigl({1\over \lambda^{\prime}}-\lambda^{\prime}\bigr),
\ \ (\lambda,\lambda^{\prime})\in T\times T.\cr}\eqno(27)$$
We use that in the variables $\lambda$, $\lambda^{\prime}$ formulas
(18), (19) take the form
$$f_y(\lambda,\lambda^{\prime},E)=f(\lambda,\lambda^{\prime},E)
\exp\bigl[{i\over 2}\sqrt{E}\bigl(\lambda\bar y+{y\over \lambda}-
\lambda^{\prime}\bar y-{y\over \lambda^{\prime}}\bigr)\bigr],\eqno(28)$$
where $(\lambda,\lambda^{\prime})\in T\times T$, $y$ is considered as
$y=y_1+iy_2$,
$$f(\lambda,\lambda^{\prime},E,t)=f(\lambda,\lambda^{\prime},E,0)
\exp\bigl[iE^{3/2}t\bigl(\lambda^3+{1\over \lambda^3}-(\lambda^{\prime})^3-
\bigl({1\over \lambda^{\prime}}\bigr)^3\bigr)\bigr],\eqno(29)$$
where  $(\lambda,\lambda^{\prime})\in T\times T$,  $t\in\R$.

The assumptions of Proposition  2 and Lemmas 1 and 2 (with (18), (19) written
 as (28), (29)) imply that
$$\eqalign{
&f(\lambda,\lambda^{\prime},E)
\exp\bigl[{i\over 2}\sqrt{E}t\bigl(\lambda\bar c+{c\over \lambda}-
\lambda^{\prime}\bar c-{c\over \lambda^{\prime}}\bigr)\bigr]=\cr
&f(\lambda,\lambda^{\prime},E)
\exp\bigl[iE^{3/2}t\bigl(\lambda^3+{1\over \lambda^3}-(\lambda^{\prime})^3-
\bigl({1\over \lambda^{\prime}}\bigr)^3\bigr)\bigr]\cr}\eqno(30)$$
for $(\lambda,\lambda^{\prime})\in T\times T$,  $t\in\R$, where $f$
is the scattering amplitude for $v(x,0)=V(x)$, $c$ is considered as
$c=c_1+ic_2$.

Property (16), identity (30) and the fact that $\lambda^3$, $\lambda^{-3}$,
$\lambda$, $\lambda^{-1}$, 1 are linear independent on each nonempty open
subset of $T$ imply (23).
\vskip 2 mm
{\it 4. Proof of Theorem 1 and final remark.}
Theorem 1 follows from  Proposition 1 and Proposition 2.

Finally, note that the result of Theorem 1 does not hold, in general,
without the assumption that $V(x)=O(e^{-\alpha |x|})$ as $|x|\to\infty$
 for some $\alpha>0$: "counter examples" to Theorem 1 with rational
bounded $V$ decaying at infinity as $O(|x|^{-2})$ are  contained (in fact)
in [G1], [G2].
As regards prototypical algebraically decaying solitons for KP1 equation,
see [FA].

\vskip 4 mm
{\bf References}
\item{[BLMP]} M.Boiti, J.J.-P.Leon, L.Martina, F.Pempinelli,
Scattering of localized solitons in the plane, Physics Letters A
{\bf 132}(8,9), 432-439 (1988)
\item{[  FA]} A.S.Fokas, M.J.Ablowitz, On the inverse scattering of the
time-dependent Schr\"odinger equation and the associated Kadomtsev-
Petviashvili (I) equation, Studies in Appl. Math. {\bf 69}, 211-228 (1983)
\item{[  FS]} A.S.Fokas, P.M.Santini, Coherent structures in multidimensions,
Phys.Rev.Lett. {\bf 63}, 1329-1333 (1989)
\item{[  G1]} P.G.Grinevich, Rational solitons of the Veselov-Novikov
equation - Two-dimensional potentials that are reflectionless for fixed
energy, Teoret. i Mat. Fiz. {\bf 69}(2), 307-310 (1986) (in Russian);
English translation: Theoret. and Math.Phys. {\bf 69}, 1170-1172 (1986).
\item{[  G2]} P.G.Grinevich, The scattering transform for the
two-dimensional Schr\"odinger operator with a potential that decreases at
infinity at fixed nonzero energy, Uspekhi Mat.Nauk {\bf 55}(6), 3-70 (2000)
(in Russian); English translation: Russian Math. Surveys {\bf 55}(6),
1015-1083 (2000)
\item{[  GN]} P.G.Grinevich, R.G.Novikov, Transparent potentials at fixed
energy in dimension two. Fixed-energy dispersion relations for the fast
decaying potentials. Cmmun.Math.Phys. {\bf 174}, 400-446 (1995).
\item{[  KN]} A.V.Kazeykina, R.G.Novikov, Absence of exponentially
localized solitons for the

\item{      } Novikov-Veselov equation at negative energy,
(in preparation)
\item{[   M]} S.V.Manakov, The inverse scattering method and two-dimensional
evolution equations, Uspekhi Mat.Nauk {\bf 31}(5), 245-246 (1976) (in Russian)
\item{[  N1]} R.G.Novikov, Construction of a two-dimensional Schr\"odinger
operator with a given scattering amplitude at fixed energy, Teoret. i Mat Fiz
{\bf 66}(2), 234-240 (1986) (in Russian); English translation: Theoret. and
Math.Phys. {\bf 66}, 154-158 (1986).
\item{[  N2]} R.G.Novikov, Reconstruction of a two-dimensional Schr\"odinger
operator from the scattering amplitude at fixed energy, Funkt.Anal. i Pril.
{\bf 20}(3), 90-91 (1986) (in Russian); English translation: Funct.Anal. and
Appl. {\bf 20}, 246-248 (1986).
\item{[  N3]} R.G.Novikov, Inverse scattering problem for the two-dimensional
Schr\"odinger equation at fixed energy and nonlinear equations, PhD Thesis,
Moscow State University 1989 (in Russian)
\item{[  N4]} R.G.Novikov, The inverse scattering problem on a fixed energy
level for the two-dimensional Schr\"odinger operator, J.Funct.Anal. {\bf 103},
 409-463 (1992).
\item{[NMPZ]} S.Novikov, S.V.Manakov, L.P.Pitaevskii, V.Z.Zakharov,
Theory of solitons: the inverse scattering method, Springer, 1984.
\item{[ NV1]} S.P.Novikov, A.P.Veselov, Finite-zone, two-dimensional,
potential Schr\"odinger operators. Explicit formula and evolution equations.
Dokl.Akad.Nauk. SSSR {\bf 279}, 20-24 (1984) (in Russian); English
translation: Sov.Math.Dokl. {\bf 30}, 588-591 (1984).
\item{[ NV2]} S.P.Novikov, A.P.Veselov, Finite-zone, two-dimensional
Schr\"odinger operators. Potential operators. Dokl.Akad.Nauk. SSSR {\bf 279},
784-788 (1984) (in Russian); English translation: Sov.Math.Dokl. {\bf 30},
705-708 (1984).
\item{[  ZS]} V.E.Zakharov, E.I.Shulman, Integrability of nonlinear systems
and perturbation theory / / What is integrability? Berlin: Springer-Verlag,
185-250 (1991)

\end